\documentclass[12pt]{article}
\usepackage{amssymb,amsmath,theorem}

\let\eps=\varepsilon \let\kappa=\varkappa                  
\newcommand\abs[1]{\left|#1\right|}                        
\newcommand\norm[1]{\left\|#1\right\|}                     
\newcommand\qed{\ifhmode\unskip\nobreak\fi\quad            
   \ifmmode\square\else\hbox{$\square$}\fi}                
\newcommand\pin{\kern.0833em}                              

\theorempreskipamount=\medskipamount \theorempostskipamount=\medskipamount

 \newcommand\proof[1]{{\it Proof#1}\,}
 \newtheorem{theorem}{Theorem}
 \newtheorem{lemma}[theorem]{Lemma}
 
 \newtheorem{remark}[theorem]{Remark}

\begin{document}

 \begin{center}
 {\large\bfseries
 \uppercase{New direct proofs of variational principles\\[1pt] for $\boldsymbol{t}$-entropy,
 spectral radius of weighted\\[1pt] shift operators, and entropy statistic\\[1pt] theorem\\}}

 \medskip

 {\large\mdseries\itshape\bigskip V.\,I.\ Bakhtin $($bakhtin@tut.by$)$$^\dag$, \
 A.\,V.\ Lebedev $($lebedev@bsu.by$)$$^\ddag$}

 \bigskip

 $^\dag$Belarusian State University / John Paul II Catholic University of Lublin, Poland;\\
 $^\ddag$\,Belarusian State University / University of Bialystok, Poland
 \end{center}

\vspace{-18pt}

\renewcommand{\abstractname}{}

\begin{abstract} \noindent
New direct proofs of variational principles for $t$-entropy, spectral radius of weighted shift
operators, and entropy statistic theorem are given. The equivalence of these statements is
obtained.
\end{abstract}

\medbreak


\smallbreak


\bigskip\medskip


As is known $t$-entropy serves as a key ingredient in evaluation of spectral radii of weighted
shift and transfer operators by means of the corresponding variational principles (see, for
example, \cite{ABL1}, \cite{ABL2}, \cite{Bakhtin2} and the sources cited therein). These
variational principles in fact show that the spectral potential (i.e., logarithm of the spectral
radius) is the Legendre transform of $t$-entropy (see \eqref{,,5} below). One of the principal
elements in all the proofs of these variational principles is the so called `entropy statistic
theorem', which proof in turn is far from being elementary (see, for example, \cite{ABL1},
\cite{Bakhtin2}). As is shown in \cite[Propositions~8.4, 8.5, 8.6]{ABL1}, $t$-entropy itself is
concave and upper semicontinuous with respect to the $^*$-weak topology (the proof here again is
nonelementary). On this base by Fenchel---Legendre---Moreau duality theorem one can establish the
corresponding variational principle for $t$-entropy, namely that $t$-entropy is the Legendre dual
to spectral potential (though this fact about $t$-entropy has not been emphasized in the mentioned
sources). To summarize the aforementioned argument, the next chain of statements has been proven:
`entropy statistic theorem' $\Rightarrow$ `variational principle for spectral potential'
$\Rightarrow$ `variational principle for $t$-entropy', where each step is rather nontrivial. In
this article we give new direct proofs for the inverted chain: `variational principle for
$t$-entropy' $\Rightarrow$ `variational principle for spectral potential' $\Rightarrow$ `entropy
statistic theorem'. In particular, it means that these statements are equivalent: `variational
principle for $t$-entropy' $\Leftrightarrow$ `variational principle for spectral potential'
$\Leftrightarrow$ `entropy statistic theorem'.

The proofs presented are shorter, more transparent and give new insight into the nature of the
objects under study (from our point of view). As a byproduct we automatically obtain concaveness
and upper semicontinuity of $t$-entropy (Remark~\ref{..8}) and also write out explicitly an
approximating sequence in the variational principle for $t$-entropy (see \eqref{,,17}).

\section{The objects of study. Variational principle for\\ spectral potential}\label{1..}

Let $(X,\mathfrak A)$ be a measurable space with a positive measure $m$ and $\alpha\!:X\to X$ be a
measurable mapping. This mapping generates a linear shift operator $A$ that maps any measurable
function $f$ to $Af = f\circ\alpha$. We will assume that this operator is bounded in the space of
integrable functions $L^1(X,m)$. The boundedness condition for $A$ is equivalent to existence of a
positive constant~$C$ such that for every measurable set $G\in \mathfrak A$ the inequality
$m(\alpha^{-1}(G))\le Cm(G)$ takes place, and in this case one has $\|A\|\le C$. For every function
$\varphi\in L^\infty(X,m)$ we define a weighted shift operator $A_\varphi$ in $L^1(X,m)$ acting
according to the formula
\begin{equation}\label{,,1}
 \bigl[A_\varphi f\bigr](x) = e^{\varphi(x)}f(\alpha(x)).
\end{equation}

 \medskip\noindent
Let $\lambda(\varphi)$ be the logarithm of its spectral radius:
\begin{equation}\label{,,2}
 \lambda(\varphi) :=\lim_{n\to\infty}\frac{1}{n}\ln\|A_\varphi^n\|,\qquad \varphi\in L^\infty(X,m).
\end{equation}
The functional $\lambda(\varphi)$ will be called the \emph{spectral potential} of operator $A$.

Recall that a linear functional on $L^\infty(X,m)$ is called positive if its values on nonnegative
functions are nonnegative, and it is normalized if its value on the unit function is equal to one.
Let $M(X,m)$ be the set of all positive normalized linear functionals on $L^\infty(X,m)$. Clearly,
$M(X,m)$ can be identified with the set of all finitely additive probability measures on
$\sigma$-algebra $\mathfrak A$ that are absolutely continuous with respect to $m$ (i.e., they take
zero values on the sets of zero measure $m$). Therefore henceforth elements of $M(X,m)$ will be
called \emph{measures}.

A measure $\mu\in M(X,m)$ is called $\alpha$-invariant if $\mu[f\circ\alpha] = \mu[f]$ for all
functions $f\in L^\infty(X,m)$. This is equivalent to the property $\mu(\alpha^{-1}(G)) = \mu(G)$
for every measurable set $G$. The family of all $\alpha$-invariant measures from $M(X,m)$ will be
denoted by~$M_\alpha(X,m)$.

By a measurable partition of unity on $X$ here we mean any finite set $D = \{g_1,\dots,g_k\}$
consisting of nonnegative functions $g_i\in L^\infty(X,m)$ satisfying the equality $g_1+\dots+g_k =
1$. For any invariant measure $\mu\in M_\alpha(X,m)$ its $t$\emph{-entropy} $\tau(\mu)$ is given by
the formulae
\begin{gather}\label{,,3}
 \tau(\mu) = \inf_{n,D}\frac{\tau_n(\mu,D)}{n}, \\[6pt]
 \tau_n(\mu,D) = \sup_{\|f\| = 1}\sum_{g\in D}\mu[g] \ln\frac{\int_X
 g\cdot|f\circ\alpha^n|\,dm}{\mu[g]}. \label{,,4}
\end{gather}
The infimum in \eqref{,,3} is taken over all natural numbers $n$ and measurable partitions of unity
$D$ on $X$ and the supremum in \eqref{,,4} is taken over all functions $f\in L^1(X,m)$ with unit
norm. If for a certain function  $g\in D$ we have $\mu[g] = 0$ then it is assumed that the
corresponding summand in~\eqref{,,4} is zero regardless of the fraction under logarithm. And if
$\mu[g]>0$ and the expression under logarithm takes zero value then it is assumed that $\mu[g]\ln 0
=-\infty$.


It is known (see \cite{ABL2}, \cite{Bakhtin2}) that spectral potential and $t$-entropy satisfy the
so-called \emph{variational principle for spectral potential} given by the formula
\begin{equation}\label{,,5}
 \lambda(\varphi) = \max_{\mu\in M_\alpha(X,m)} \bigl(\tau(\mu) +\mu[\varphi]\bigr), \qquad
 \varphi\in L^\infty(X,m).
\end{equation}

We will give a new proof for this formula in Theorem~\ref{..10}. It will be a simple corollary of
the dual result (variational principle for $t$-entropy); namely its new direct proof is the
starting point of the next section.

\section{New proofs of variational principles for $\boldsymbol{t}$-entropy\\
and spectral potential}\label{2..}

\begin{theorem}[variational principle for $\boldsymbol{t}$-entropy] \label{..1}
The following equality takes place\/$:$
\begin{equation} \label{,,6}
 \tau(\mu) =\inf_{\varphi\in L^\infty(X,m)} \bigl(\lambda(\varphi) -\mu[\varphi]\bigr), \qquad
 \mu\in M_\alpha(X,m).
\end{equation}
\end{theorem}

\begin{theorem} \label{..2}
For each linear functional\/ $\mu$ on the space\/ $L^\infty(X,m)$ that does not belong to\/
$M_\alpha(X,m)$ the following equality takes place\/$:$
\begin{equation} \label{,,7}
 \inf_{\varphi\in L^\infty(X,m)} \bigl(\lambda(\varphi) -\mu[\varphi]\bigr) = -\infty.
\end{equation}
\end{theorem}

These two theorems show that it is natural to put
\begin{equation} \label{,,8}
 \tau(\mu) =-\infty \quad \text{for all} \ \
 \mu\in L^\infty(X,m)\rule{0ex}{1.7ex}^*\setminus M_\alpha(X,m).
\end{equation}
Then formulae \eqref{,,6} and \eqref{,,7} are united into one:
\begin{equation} \label{,,9}
 \tau(\mu) =\inf_{\varphi\in L^\infty(X,m)} \bigl(\lambda(\varphi) -\mu[\varphi]\bigr), \qquad
 \mu\in L^\infty(X,m)\rule{0ex}{1.7ex}^*.
\end{equation}

We start the proof of Theorem~\ref{..1} with the next lemma presenting, in particular, `Young'
inequality (inequality \eqref{,,11}) for the objects considered. Here we take Young in quotation
marks since the classical Young inequality (serving as a component of Fenchel---Legendre duality
construction) exploits convexity and Legendre dual objects, while up to now and so also in the
proof of Lemma~\ref{..3} we do not use them. In fact inequality \eqref{,,11} will become the
classical Young inequality on derivation of Theorem~\ref{..1} (cf. Remark~\ref{..9}).

\begin{lemma} \label{..3}
For any $\varphi \in L^\infty(X,m)$, $\mu\in M_\alpha(X,m)$, $n\in\mathbb N$ and $\eps >0$  there
exists a measurable partition of unity $D$ such that
\begin{equation} \label{,,10}
 \eps + \frac{\ln\|A_\varphi^n\|}{n}
 \ge \mu[\varphi] +\frac{\tau_n(\mu,D)}{n}
\end{equation}
and thus
\begin{equation} \label{,,11}
 \lambda(\varphi) \ge \mu[\varphi] +\tau(\mu).
\end{equation}
\end{lemma}

 \medskip

\proof. The proof, in fact, can be extracted from the reasoning in \cite{Bakhtin2} and we give it
for completeness of presentation.

Note first that by arbitrariness of $\eps >0$ inequality \eqref{,,11} follows from \eqref{,,10} by
taking infimum with respect to $n,\,D$.

So it is enough to verify \eqref{,,10}.

Let us introduce the notation
\begin{equation*}
 S_n\varphi = \varphi +\varphi\circ\alpha+ \,\dotsm\,+ \varphi\circ \alpha^{n-1}.
\end{equation*}
Then \eqref{,,1} implies the equality
\begin{equation*}
 A_\varphi^n f = e^{S_n\varphi} f\circ \alpha^n.
\end{equation*}

 \medskip\noindent
For arbitrary numbers $n\in\mathbb N$ and $\eps>0$ we choose a fine measurable partition of
unity~$D$ such that on the support of each function $g\in D$ the essential oscillation of function
$S_n\varphi$ does not exceed~$\eps$. This~$D$ is in fact the desired partition.

Indeed. Let $S_n\varphi(g) := \mathop{\mathrm{ess\,sup}}\pin\{\pin S_n\varphi(x)\mid g(x)\ne
0\pin\}$. Then for any $\mu\in M_\alpha(X,m)$ and integrable function $f$ concavity of logarithm
implies the inequalities
\begin{align*}
\eps+\ln\int_X e^{S_n\varphi}\abs{f\circ\alpha^n}dm &\ge \ln\sum_{g\in D}e^{S_n\varphi(g)}
\int_X g\abs{f\circ\alpha^n}dm\\[3pt]
&\ge \sum_{g\in D}\mu[g]\ln\frac{e^{S_n\varphi(g)}
\int_X g\abs{f\circ\alpha^n}dm}{\mu[g]}\\[3pt]
&\ge \mu[S_n\varphi] +\sum_{g\in D}\mu[g]\ln \frac{\int_X g\abs{f\circ\alpha^n}dm}{\mu[g]}.
\end{align*}
Passing here to supremum with respect to functions $f$ with unit norm and taking into account
\eqref{,,4}, one obtains the inequality
\begin{equation*}
 \eps+ \ln\|A_\varphi^n\|\ge \mu[S_n\varphi] +\tau_n(\mu,D) =n\mu[\varphi] +\tau_n(\mu,D),
\end{equation*}
which implies \eqref{,,10}. \qed

\medskip

Inequality \eqref{,,11} implies that  the left hand part in \eqref{,,6} does not exceed its right
hand part.

To finish the proof of Theorem~\ref{..1} we need two more lemmas.

In the next lemma notation $\lambda(\varphi,A)$ has the same meaning as $\lambda(\varphi)$ and
$\lambda(n\varphi,A^n)$ denotes logarithm of the spectral radius of operator $e^{n\varphi}\!\pin
A^n$.

\begin{lemma} \label{..4}
The following inequality takes place\/$:$
\begin{equation} \label{,,12}
 n\lambda(\varphi,A) \le \lambda(n\varphi,A^n), \qquad n\in \mathbb N.
\end{equation}
\end{lemma}

\proof. Note that for any natural $k$ one has
\begin{equation*}
 \exp\{S_{nk}\varphi\} =
 \exp\biggl\{\sum_{i=0}^{n-1}\sum_{j=0}^{k-1} \varphi\circ\alpha^{i+nj}\biggr\} =
 \prod_{i=0}^{n-1} \exp\biggl\{\sum_{j=0}^{k-1} \varphi\circ\alpha^{i+nj}\biggr\}.
\end{equation*}
Let $c =\mathop{\mathrm{ess\,sup}}|\varphi|$. Exploiting H\"older inequality in the form
\begin{equation*}
 \int_X \psi_1 \dotsm \psi_n \, d\nu \le \prod_{i=1}^{n}\left(\int_X \abs{\psi_i}^n
 d\nu\right)^{\!1/n},
\end{equation*}
where $d\nu = \abs{f\circ\alpha^{n(k+1)}} dm$, one obtains
\begin{align*}
 e^{-nc}\pin\big\|A_\varphi^{n(k+1)}f\big\| &=
 e^{-nc}\!\int_X e^{S_{n(k+1)}\varphi}\abs{f\circ\alpha^{n(k+1)}} dm \le
 \int_X e^{S_{nk}\varphi}\abs{f\circ\alpha^{n(k+1)}} dm \\[6pt]
 &\le \prod_{i=0}^{n-1} \bigg(\!\pin\int_X
 \exp\biggl\{\sum_{j=0}^{k-1} (n\varphi\circ\alpha^{nj})\circ\alpha^i\biggr\}
 \pin\big|(f\circ\alpha^{n-i})\circ\alpha^{nk}\circ\alpha^i\big|\,dm\bigg)^{\!1/n} \\[4pt]
 &= \prod_{i=0}^{n-1} \big\|A^i(e^{n\varphi}\!\pin A^n)^k(A^{n-i}f)\big\|^{1/n} \le
 \norm{A}^n\big\|(e^{n\varphi}\!\pin A^n)^k\big\| \norm{f}.
\end{align*}
Thus,
\begin{equation*}
 -nc +\ln\big\|A_\varphi^{n(k+1)}\big\| \le n\ln\|A\| +\ln\big\|(e^{n\varphi}\!\pin A^n)^k\big\|.
\end{equation*}

 \medskip\noindent
Dividing the latter inequality by $k$ and turning $k\to\infty$ one gets  \eqref{,,12}. \qed

\medskip

Let us fix a measure $\mu\in M_\alpha(X,m)$, natural number $n$ and measurable partition of unity
$D$ on $X$. For these objects there exists a sequence of functions $f_k\in L^1(X,m)$ with unit norm
on which the supremum in \eqref{,,4} is attained. One may choose a subsequence $f_{k_i}$ of this
sequence such that the following limits do exist simultaneously:
\begin{equation} \label{,,13}
 \lim_{i\to\infty} \int_X g \abs{f_{k_i}\!\circ\alpha^n} dm =: C_n(\mu,g,D), \qquad g\in D.
\end{equation}
Then by construction one has
\begin{equation} \label{,,14}
 \tau_n(\mu,D) =
 \sum_{g\in D}\mu[g] \ln\frac{C_n(\mu,g,D)}{\mu[g]}.
\end{equation}

\vspace{-1pt}

\begin{lemma} \label{..5}
If\/ $\tau_n(\mu,D)>-\infty$ then
\begin{equation} \label{,,15}
 \sup_{\|f\|=1} \sum_{\substack{g\in D,\\[1pt] \mu[g]>0}}
 \mu[g]\pin \frac{\int_X g \abs{f\circ\alpha^n} dm}{C_n(\mu,g,D)} \,=\, 1.
\end{equation}
\end{lemma}

\proof. Finiteness of $\tau_n(\mu,D)$ and \eqref{,,14} imply that $C_n(\mu,g,D)>0$ whenever
$\mu[g]>0$. For each $f\in L^1(X,m)$ with $\|f\| =1$ let us consider the function
\begin{equation*}
 \eta(t) =\sum_{\substack{g\in D,\\[1pt] \mu[g]>0}} \mu[g]\ln
 \frac{(1-t)C_n(\mu,g,D) +t\int_X g\abs{f\circ \alpha^n} dm}{\mu[g]}, \qquad t\in [0,1].
\end{equation*}
By definition of the numbers $C_n(\mu,g,D)$ this function attains its maximal value equal to
$\tau_n(\mu,D)$ at $t=0$. Therefore its derivative at $t=0$
\begin{equation*}
 \frac{d\eta(t)}{dt}\bigg|_{t=0} =\sum_{\substack{g\in D,\\[1pt] \mu[g]>0}}
 \mu[g]\pin \frac{\int_X g\abs{f\circ\alpha^n} dm -C_n(\mu,g,D)}{C_n(\mu,g,D)} =
 \sum_{\substack{g\in D,\\[1pt] \mu[g]>0}}
 \mu[g]\pin \frac{\int_X g\abs{f\circ\alpha^n} dm}{C_n(\mu,g,D)} -1
\end{equation*}
is nonpositive and so the left hand part in \eqref{,,15} does not exceed its right hand part.

The equality in \eqref{,,15} is attained on the sequence of functions $f_{k_i}$ from \eqref{,,13}.
\qed

\medskip

Now we can finish the proof of Theorem~\ref{..1}.

Let us fix an arbitrary measure $\mu\in M_\alpha(X,m)$, natural number $n$ and measurable partition
of unity $D$ on $X$.

Suppose at first that there exists a function $g\in D$ satisfying the inequality  $\mu[g] >0$ and
equality $\int_X g\abs{f\circ\alpha^n} dm =0$ for all $f\in L^1(X,m)$. Then by definition one has
$\tau_n(\mu,D) =-\infty$ and therefore $\tau(\mu) =-\infty$. Thus in this case equality \eqref{,,6}
takes the form
\begin{equation} \label{,,16}
 -\infty =\inf_{\varphi\in L^\infty(X,m)} \bigl(\lambda(\varphi) -\mu[\varphi]\bigr).
\end{equation}
Let us verify it.

Consider the family of functions $\varphi_t =tg/n$ where $t\in\mathbb R$. Inequalities $0\le g\le
1$ and the Lagrange theorem imply that
\begin{equation*}
 e^{n\varphi_t} =e^{tg} \le 1+e^t tg.
\end{equation*}
Therefore for each function $f\in L^1(X,m)$ with unit norm one has
\begin{equation*}
 \|e^{n\varphi_t}\!\pin A^nf\| =\int_X e^{n\varphi_t}|f\circ\alpha^n|\,dm \le
 \int_X (1+e^ttg)\pin |f\circ\alpha^n|\,dm =\int_X |f\circ\alpha^n|\,dm \le \|A^n\|.
\end{equation*}
Thus $\|e^{n\varphi_t}\!\pin A^n\| \le \|A^n\|$. Applying Lemma~\ref{..4}, we obtain the following
estimate
\begin{equation*}
 n\lambda(\varphi_t) =n\lambda(\varphi_t,A) \le \lambda(n\varphi_t,A^n) \le
 \ln \|e^{n\varphi_t}\!\pin A^n\| \le \ln \|A^n\|.
\end{equation*}
On the other hand,
\begin{equation*}
 \mu[\varphi_t] =\mu[tg/n] = t\mu[g]/n \to +\infty \quad\text{as}\ \ t\to +\infty.
\end{equation*}
And therefore $\lambda(\varphi_t) -\mu[\varphi_t] \to -\infty$ when $t\to +\infty$. So equality
\eqref{,,16} is verified.

\smallskip

It remains to consider the situation when for each function $g\in D$ satisfying the condition
$\mu[g]
>0$ there exists a function $f_g\in L^1(X,m)$ such that $\int_X g\abs{f_g\circ\alpha^n} dm
>0$. Taking the function  $f :=\sum_g |f_g|$
one obtains that
\begin{equation*}
 \int_X g\abs{f\circ\alpha^n} dm >0 \quad\text{as soon as }\ \ \mu[g]>0.
\end{equation*}
Therefore  $\tau_n(\mu,D) >-\infty$. Note also that finiteness of  $\tau_n(\mu,D)$ along with
\eqref{,,14} implies that the condition $\mu[g]>0$ automatically implies the inequality
$C_n(\mu,g,D)>0$.

Now let us define the family of functions
\begin{equation}\label{,,17}
 \varphi_\eps  :=\, \frac{1}{n}\pin\ln\Bigg\{
 \sum_{\mu[g]>0} \frac{\mu[g]}{C_n(\mu,g,D)}\pin g +
 \sum_{\mu[g]=0} \eps g\Bigg\}, \qquad \eps>0.
\end{equation}
For any integrable function $f$ with unit norm one has
\begin{align*}
 \|e^{n\varphi_\eps}\!\pin A^n f\| &=\int_X e^{n\varphi_\eps}|f\circ\alpha^n|\, dm \\[6pt]
 &= \int_X \sum_{\mu[g]>0} \frac{\mu[g]}{C_n(\mu,g,D)}\pin g\abs{f\circ\alpha^n} dm +
 \int_X \sum_{\mu[g]=0} \eps g\abs{f\circ\alpha^n} dm \\[6pt]
 &=\sum_{\mu[g]>0} \mu[g]\pin \frac{\int_X g \abs{f\circ\alpha^n} dm}{C_n(\mu,g,D)} +
 \eps\int_X \sum_{\mu[g]=0} g\abs{f\circ\alpha^n} dm \le 1 +\eps\|A^n\|
\end{align*}
(where in the final inequality we exploited Lemma~\ref{..5}). This along with Lemma~\ref{..4}
implies the estimate
\begin{equation} \label{,,18}
 n\lambda(\varphi_\eps) \le \lambda(n\varphi_\eps,A^n) \le
 \ln \|e^{n\varphi_\eps}\!\pin A^n\| \le \ln(1+\eps\|A^n\|) \le \eps\|A^n\|.
\end{equation}

On the other hand applying concaveness of logarithm and \eqref{,,14} one obtains
\begin{equation} \label{,,19}
\begin{aligned}
 \mu[n\varphi_\eps] &\,=\,\mu\Bigg[\ln\Bigg\{\sum_{\mu[g]>0} \frac{\mu[g]}{C_n(\mu,g,D)}\pin g +
 \sum_{\mu[g]=0} \eps g\Bigg\}\Bigg] \\[6pt]
 &\,\ge\, \mu\Bigg[\sum_{\mu[g]>0} g\ln\frac{\mu[g]}{C_n(\mu,g,D)} +
 \sum_{\mu[g]=0} g\ln\eps\Bigg] =\,-\tau_n(\mu,D).
\end{aligned}
\end{equation}
Combining \eqref{,,19} and \eqref{,,18} we get
\begin{equation*}
 \frac{\tau_n(\mu,D)}{n} \ge -\mu[\varphi_\eps] \ge
 -\mu[\varphi_\eps] + \left(\lambda(\varphi_\eps)
 -\frac{\eps\|A^n\|}{n}\right)
 \end{equation*}
and therefore
 \begin{equation*}
 \frac{\tau_n(\mu,D)}{n} +\frac{\eps\|A^n\|}{n} \ge
 \lambda(\varphi_\eps)
 -\mu[\varphi_\eps] \ge
 \inf_{\varphi\in L^\infty(X,m)}\big(\lambda(\varphi) -\mu[\varphi]\big).
\end{equation*}
This inequality along with arbitrariness of $\eps$, $n$, $D$ and definition \eqref{,,3} of
$\tau(\mu)$ implies the inequality
\begin{equation*}
\tau(\mu) \ge \inf_{\varphi\in L^\infty(X,m)}\big(\lambda(\varphi) -\mu[\varphi]\big).
\end{equation*}

 \medskip\noindent
Together with Young inequality \eqref{,,11} this proves  \eqref{,,6}. \qed

\medskip

Now let us prove Theorem~\ref{..2}. For this purpose we need two more lemmas describing the
principal properties of the functional $\lambda(\varphi)$.

\begin{lemma} \label{..6}
Spectral potential\/ $\lambda(\varphi)$ possesses the following properties\/$:$

\smallskip

a\/$)$ if\/ $\varphi\ge\psi$, then\/ $\lambda(\varphi)\ge\lambda(\psi)$ \,\emph{(monotonicity);}

\smallskip

b\/$)$ \,$\lambda(\varphi+t) =\lambda(\varphi)+t$ for all\/ $t\in\mathbb R$ \,\emph{(additive
homogeneity);}

\smallskip

c\/$)$ \,$|\lambda(\varphi) -\lambda(\psi)|\le \mathop{\mathrm{ess\,sup}} |\varphi-\psi|$
\,\emph{(Lipschitz condition);}

\smallskip

d\/$)$ \,$\lambda((1-t)\varphi +t\psi)\le (1-t)\lambda(\varphi) +t\lambda(\psi)$ for\/ $t\in [0,1]$
\,\emph{(convexity);}

\smallskip

e\/$)$ \,$\lambda(\varphi +\psi\circ\alpha) =\lambda(\varphi+\psi)$ \,\emph{(strong
$\alpha$-invariance)}.

\end{lemma}

This lemma is proven in \cite{ABL1}, \cite{Bakhtin1}.

\begin{lemma} \label{..7}
If a linear functional\/ $\mu$ on\/ $L^\infty(X,m)$ possesses the property
\begin{equation} \label{,,20}
 \inf_{\varphi\in L^\infty(X,m)} \big(\lambda(\varphi) -\mu[\varphi]\big) >-\infty,
\end{equation}
then\/ $\mu\in M_\alpha(X,m)$. In particular, this is true for every subgradient of the function\/
$\lambda(\varphi)$.
\end{lemma}

\proof. Observe first that any functional $\mu$ satisfying \eqref{,,20} is necessarily positive.
Indeed, if $\mu[\varphi]<0$ for some function $\varphi\ge 0$ then by monotonicity of
$\lambda(\varphi)$ (Lemma~\ref{..6}, a)\,) one has $\lambda(-t\varphi)\le \lambda(0)$ and thus
\begin{equation*}
 \lambda(-t\varphi) -\mu[-t\varphi] \le \lambda(0) +t\mu[\varphi] \to -\infty \quad\text{as}\ \
 t\to +\infty,
\end{equation*}
which contradicts  \eqref{,,20}.

Note also that this functional should have norm equal to 1. Indeed, by additive homogeneity of
$\lambda(\varphi)$ (Lemma~\ref{..6}, b)\,) we have $\lambda(t) =\lambda(0) +t$ and therefore
\begin{equation*}
 \lambda(t) -\mu[t] =\lambda(0) +t(1-\mu[1]).
\end{equation*}
This function is bounded from below for all $t\in\mathbb R$ only if $\mu[1] =1$. Which by
positivity of $\mu$ means that $\Vert \mu \Vert =1$.

Observe now that strong $\alpha$-invariance of $\lambda(\varphi)$ (Lemma~\ref{..6}, e)\,) implies
the equality $\lambda(t\varphi\circ\alpha -t\varphi) = \lambda(0)$. Thus,
\begin{equation*}
 \lambda(t\varphi\circ\alpha -t\varphi) -\mu[t\varphi\circ\alpha -t\varphi] =
 \lambda(0) -t(\mu[\varphi\circ\alpha] -\mu[\varphi]).
\end{equation*}
This function is bounded from below for all $t\in\mathbb R$ only if $\mu[\varphi\circ\alpha]
=\mu[\varphi]$ which means $\alpha$-invariance of $\mu$.

\smallskip

Thus condition \eqref{,,20} can be satisfied only when $\mu\in M_\alpha(X,m)$.

\smallskip

Finally, we observe that any subgradient  $\mu$ of $\lambda(\varphi)$ satisfies \eqref{,,20}.
Indeed, let $\mu$ be a subgradient of $\lambda(\varphi)$ at a point $\varphi_0$. Then
\begin{equation*}
 \lambda(\varphi) -\lambda(\varphi_0) \ge \mu[\varphi-\varphi_0] \quad\text{for all}\ \
 \varphi\in L^\infty(X,m),
\end{equation*}
that is
\begin{equation*}
 \lambda(\varphi) -\mu[\varphi] \ge \lambda(\varphi_0) -\mu[\varphi_0] \quad\text{for all}\ \
 \varphi\in L^\infty(X,m),
\end{equation*}

 \medskip\noindent
which implies \eqref{,,20}. \qed

\medskip

Clearly, Theorem~\ref{..2} is a straightforward corollary of Lemma~\ref{..7}.

\begin{remark} \label{..8}
{\rm By Lemma~\ref{..6} the function $\lambda(\varphi)$ is convex and continuous.
Theorems~\ref{..1} and \ref{..2} in essence state that the functional $-\tau(\mu)$ is the Legendre
transform of $\lambda(\varphi)$. This automatically implies that $t$-entropy $\tau(\mu)$ is concave
and upper semicontinuous (in the $^*$-weak topology) on the dual space to $L^\infty(X,m)$. In
\cite{Bakhtin2} concaveness and upper semicontinuity of $t$-entropy were proven independently and
in an essentially more complicated way.}
\end{remark}

\begin{remark} \label{..9}
{\rm Now we see that inequality \eqref{,,11} is indeed the classical Young inequality for the
objects considered.}
\end{remark}

Finally we observe that variational principle \eqref{,,5} can be easily derived from
Theorem~\ref{..1} and Lemmas~\ref{..6}, \ref{..7}.

\begin{theorem}[variational principle for spectral potential] \label{..10}
For each function $\varphi\in L^\infty(X,m)$ the following equality takes place\/$:$
\begin{equation} \label{,,21}
\lambda(\varphi) = \max_{\mu\in M_\alpha(X,m)} \bigl(\tau(\mu) +\mu[\varphi]\bigr).
\end{equation}
\end{theorem}

\proof. By Lemma~\ref{..6} the function $\lambda(\varphi)$ is convex and continuous. Thus at each
point $\varphi_0$ there exists at least one subgradient $\mu$ for $\lambda(\varphi)$. By
Lemma~\ref{..7} this subgradient belongs to $M_\alpha(X,m)$. By Theorem~\ref{..1} and definition of
a subgradient we have
\begin{equation*}
  \tau(\mu) =\inf_{\varphi\in L^\infty(X,m)} \big(\lambda(\varphi) -\mu[\varphi]\big) =
  \lambda(\varphi_0) -\mu[\varphi_0].
\end{equation*}
Therefore $\lambda(\varphi_0) =\tau(\mu) +\mu(\varphi_0)$. Combining this equality with Young
inequality~\eqref{,,11} one obtains \eqref{,,21}. \qed

\section{Variational principles in the spaces of bounded\\ measurable functions}\label{3..}

In the above part of the article we considered the spectral potential $\lambda(\varphi)$ as a
function on the space $L^\infty(X,m)$ and $t$-entropy $\tau(\mu)$ as a function on the dual space
$L^\infty(X,m)\rule{0ex}{1.7ex}^*$. But in fact the spectral potential $\lambda(\varphi)$ can
readily be defined in the same way (i.e., by formula~\eqref{,,2}) on the space of bounded
measurable functions $B(X,\mathfrak A)$ (i.e., not on the factor space $L^\infty(X,m)$ of
$B(X,\mathfrak A)$) and $\tau(\mu)$ can be defined by the same formulae on the dual space
$\displaystyle B(X,\mathfrak A)^*$. In this case all the above presented results and their proofs
remain in force. In particular, in the situation mentioned equality~\eqref{,,9} takes form
\begin{equation} \label{,,22}
 \tau(\mu) =\inf_{\varphi\in B(X,\mathfrak A)} \bigl(\lambda(\varphi) -\mu[\varphi]\bigr), \qquad
 \mu\in B(X,\mathfrak A)^*,
\end{equation}
and variational principle~\eqref{,,21} looks as
\begin{equation} \label{,,23}
 \lambda(\varphi) =\max_{\mu\in M_\alpha(X,\mathfrak A)}\big(\tau(\mu) +\mu[\varphi]\big), \qquad
 \varphi\in B(X,\mathfrak A),
\end{equation}
where $M_\alpha(X,\mathfrak A)$ denotes the set of all $\alpha$-invariant positive normalized
linear functionals on $B(X,\mathfrak A)$.

Since the natural factor mapping from $\displaystyle B(X,\mathfrak A)$ onto $L^\infty(X,m)$ is
linear and bounded it induces a continuous embedding of
 $L^\infty(X,m)\rule{0ex}{1.7ex}^*$ in $\displaystyle
B(X,\mathfrak A)^*$. Therefore from the formal point of view  equality \eqref{,,22} is stronger
than \eqref{,,9}, and equality \eqref{,,23} is weaker than \eqref{,,21} (here we bear in mind Young
inequality \eqref{,,11} which, of course, is also true for $\varphi \in B(X,\mathfrak A)$ and $
\mu\in M_\alpha(X,\mathfrak A)$).

\smallskip

The next lemma shows that in fact all the corresponding variational principles are equivalent.

\begin{lemma} \label{..11}
If\/ $\displaystyle \mu\in B(X,\mathfrak A)^*\setminus L^\infty(X,m)\rule{0ex}{1.7ex}^*$ then
\begin{equation} \label{,,24}
 \tau(\mu) =-\infty \quad\mathrm{and}\quad
 \inf_{\varphi\in B(X,\mathfrak A)} \bigl(\lambda(\varphi)-\mu[\varphi]\bigr) =-\infty.
\end{equation}
\end{lemma}

\proof. The assumption $\displaystyle \mu\in B(X,\mathfrak A)^*\setminus
L^\infty(X,m)\rule{0ex}{1.7ex}^*$ implies that
 there exists a function $\varphi\in B(X,\mathfrak A)$
 which is equal to zero
$m$-almost everywhere and therewith  $\mu[\varphi] \ne 0$. For this function we have
\begin{equation*}
 \inf_{t\in\mathbb R} \big(\lambda(t\varphi) -\mu[t\varphi]\big) =
 \inf_{t\in\mathbb R} \big(\lambda(0) -t\mu[\varphi]\big) =-\infty
\end{equation*}
and therefore
\begin{equation*}
 \inf_{\varphi\in B(X,\mathfrak A)} \bigl(\lambda(\varphi)-\mu[\varphi]\bigr) =-\infty.
\end{equation*}

 \medskip\noindent
Thus to finish the proof it suffice to verify that $\tau(\mu) = -\infty$ for the mentioned $\mu$.

\smallskip

If $\displaystyle \mu\in B(X,\mathfrak A)^*\setminus M_\alpha(X,\mathfrak A)$ then by definition
$\tau(\mu) =-\infty$ (by analogy with \eqref{,,8}). And if $\mu\in M_\alpha(X,\mathfrak A)
\setminus L^\infty(X,m)\rule{0ex}{1.7ex}^*$ then there exists a measurable set $G\in\mathfrak A$
such that $m(G) =0$ and therewith $\mu(G) >0$. Let us take a partition of unity $D=\{g,\,1-g\}$,
where $g$ is the characteristic function of this set $G$. In view of  \eqref{,,4} one has
$\tau_n(\mu,D) =-\infty$ and therefore $\tau(\mu) =-\infty$. \qed

\section{A new proof of entropy statistic theorem}\label{4..}

As it was noted, `entropy statistic theorem' served as a key instrument in the initial proof of
variational principle for spectral potential in \cite{Bakhtin2}. It appears that its proof can also
be essentially simplified by means of the variational principle for $t$-entropy. In particular,
this means that we obtain an equivalence:

`variational principle for $t$-entropy' $\Leftrightarrow$ `variational principle for spectral
potential' $\Leftrightarrow$ `entropy statistic theorem'.

It will be covenient for us to present a variant of `entropy statistic theorem' in the spaces
$B(X,\mathfrak A)$ and $\displaystyle B(X,\mathfrak A)^*$.

Let us take an arbitrary point $x\in X$ and a natural number $n$. The corresponding \emph{empirical
measure} $\delta_{x,n}$ is concentrated at points $x$, $\alpha(x)$, \dots, $\alpha^{n-1}(x)$ and
its value at each point is equal to $1/n$. For any function $f\!:X\to\mathbb R$ its integral with
respect to this measure equals to the `empirical mean'
\begin{equation*}
 \delta_{x,n}[f] =\frac{f(x)+f(\alpha(x))+\,\dotsm\,+f(\alpha^{n-1}(x))}{n} =
 \frac{1}{n}\pin S_nf(x).
\end{equation*}
Clearly, $\delta_{x,n}$ can be considered as a linear functional on the space $B(X,\mathfrak A)$.

Let $\mu$ be a  linear functional  on $B(X,\mathfrak A)$ and $O(\mu)$ be its certain neighborhood
in the $^*$-weak topology. By $X_n(O(\mu))$ we denote the set
\begin{equation*}
 X_n(O(\mu)) :=\{\pin x\in X\mid \delta_{x,n}\in O(\mu)\pin\}.
\end{equation*}

`Entropy statistic theorem' is formulated as follows.

\begin{theorem} \label{..12}
For any linear functional\/ $\mu$ on\/ $B(X,\mathfrak A)$ and any number\/ $\eps>0$ there exists a
neighborhood\/ $O(\mu)$ in the\/ $^*$-weak topology and a large number\/ $C(\eps,\mu)$ such that
for all natural numbers\/ $n$ and functions\/ $f\in L^1(X,m)$ the next estimate holds\/$:$
\begin{equation} \label{,,25}
 \intop_{X_n(O(\mu))} \hspace{-0.4em} |f\circ\alpha^n|\,dm \,\le\,
 C(\eps,\mu)\pin e^{n(\tau(\mu)+\eps)} \intop_{X} |f|\,dm.
\end{equation}
\end{theorem}

If $\tau(\mu) =-\infty$ then the number $\tau(\mu)+\eps$ in \eqref{,,25} should be replaced by
$-1/\eps$.

\medskip

\proof. By the variational principle for $t$-entropy in the form of \eqref{,,22} one can take a
function $\varphi\in B(X,\mathfrak A)$ such that
\begin{equation*}
 \lambda(\varphi) -\mu[\varphi] < \tau(\mu)+\eps/2
\end{equation*}
(or $\lambda(\varphi) -\mu[\varphi] <-1/\eps -\eps/2$ in the case when $\tau(\mu) =-\infty)$. Let
us set
\begin{equation*}
 O(\mu) :=\{\pin \delta\in B(X,\mathfrak A)^*\mid
 \lambda(\varphi) -\delta[\varphi] <\tau(\mu) +\eps/2\pin\}.
\end{equation*}
Then for each point $x\in X_n(O(\mu))$ the following relations are true
\begin{equation*}
 S_n\varphi(x) =n\delta_{x,n}[\varphi] >n(\lambda(\varphi) -\tau(\mu) -\eps/2).
\end{equation*}
Take a constant $C(\eps,\mu)$ so large that
\begin{equation*}
 \|A_\varphi^n\| \le C(\eps,\mu)\pin e^{n(\lambda(\varphi) +\eps/2)}, \qquad n\in\mathbb N.
\end{equation*}
Now \eqref{,,25} follows from the calculation
\begin{align*}
 C(\eps,\mu)\pin e^{n(\lambda(\varphi)+\eps/2)} \|f\| \,&\ge\, \|A_\varphi^n f\| \,=\pin
 \intop_X e^{S_n\varphi}\abs{f\circ\alpha^n} dm \,\ge
 \intop_{X_n(O(\mu))} \hspace{-0.4em} e^{S_n\varphi}\abs{f\circ\alpha^n} dm \\[6pt]
 &\ge\, e^{n(\lambda(\varphi) -\tau(\mu) -\eps/2)} \hspace{-0.4em}
 \intop_{X_n(O(\mu))} \hspace{-0.4em} |f\circ\alpha^n|\, dm. \qed
\end{align*}



\end{document}